\documentclass{amsart}

\usepackage{amsmath,amssymb,amscd,amsthm,mathdots,indentfirst}
\usepackage{amsfonts}
\usepackage[mathscr]{eucal}
\usepackage{mathpazo}
\usepackage{ifpdf,ifdraft}

\ifpdf
\usepackage[pdftex,bookmarks,colorlinks,breaklinks]{hyperref}  
\hypersetup{linkcolor=red,citecolor=blue,filecolor=magenta,urlcolor=blue} 
\hypersetup{linkcolor=black,citecolor=black,filecolor=black,urlcolor=black} 
\fi

\swapnumbers
\newtheorem{prop}[equation]{Proposition}
\newtheorem{thm}[equation]{Theorem}
\newtheorem{cor}[equation]{Corollary}
\newtheorem{lem}[equation]{Lemma}

\theoremstyle{definition}
\newtheorem{defn}[equation]{Definition}

\newtheorem{remark}[equation]{Remark}

\def\map{\operatorname{map}}

\def\Mor{\mathrm{Mor}}

\def\Obj{\mathrm{Obj}}

\def\vp{\varphi}
\def\al{\alpha}
\def\be{\beta}
\def\ga{\gamma}

\def\de{\delta}
\def\De{\Delta}
\def\si{\sigma}

\def\A{\ensuremath{\mathcal{A}}}
\def\B{\ensuremath{\mathcal{B}}}
\def\C{\ensuremath{\mathcal{C}}}

\def\E{\ensuremath{\mathcal{E}}}

\def\G{\ensuremath{\mathcal{G}}}

\def\M{\ensuremath{\mathcal{M}}}

\def\O{\ensuremath{\mathcal{O}}}
\def\SS{\ensuremath{\mathcal{S}}}

\def\RR{\ensuremath{\mathbb{R}}}

\def\opEND{\operatorname{\O End}}

\def\Hom{\mathrm{Hom}}

\DeclareMathOperator{\id}{id}
\def\Id{\textrm{Id}}
\def\op{\ensuremath{\mathrm{op}}}

\def\Lan{\operatornamewithlimits{LKan}}

\def\sets{\mathbf{Sets}}

\def\Tot{\operatorname{Tot}}

\newcommand{\xto}[1]{\xrightarrow{#1}}

\DeclareMathOperator{\End}{End}

\DeclareMathOperator{\pro}{pro}

\newcommand{\comp}{\mathbin{\raise0.2ex\hbox{\tiny$\diamond$}}}

\ifdraft{

}

\ifpdf
\else
\def\texorpdfstring#1#2{#1}
\fi

\input xy
\xyoption{all}

\def\smashop#1_#2{%
\displaystyle{#1_{%
\hbox to 0pt{\hss$\scriptstyle{#2}$\hss}}\;}}

\def\smashopsup#1^#2{%
\displaystyle{#1^{%
\hbox to 0pt{\hss$\scriptstyle{#2}$\hss}}\;}}

\def\De{\Delta}

\DeclareMathOperator{\TW}{TW}

\date{December 8, 2008}
\title{A simplicial $A_\infty$-operad acting on $R$-resolutions}
\author{Tilman Bauer}
\address{Faculteit der exacte wetenschappen, Vrije Universiteit Amsterdam, de Boelelaan 1081A, 1081HV Amsterdam, The Netherlands}
\email{tilman@few.vu.nl}

\author{Assaf Libman}
\address{
	Department of Mathematical Sciences, 
	King's College,
	University of Aberdeen,
   	Aberdeen AB24 3UE , Scotland, U.K.}
\email{assaf@maths.abdn.ac.uk}


\begin{document}

\pagestyle{plain}

\numberwithin{equation}{section}
\renewcommand{\theequation}{\thesection.\arabic{equation}}
\renewcommand{\thethmain}{\Alph{thmain}}

\begin{abstract}
We construct a combinatorial model of an $A_\infty$-operad which acts simplicially on the cobar resolution (not just its total space) of a simplicial set with respect to a ring $R$.
\end{abstract}

\maketitle

\section{Introduction}
Let $K$ be a commutative ring and $X$ a simplicial set (henceforth, simply called ``space''). Bousfield and Kan \cite{yellowmonster} defined the \emph{$K$-resolution} of $X$ to be (a variant of) the cosimplicial space
\[
R_K(X)^n = \underset{n+1\text{ times}}{\underbrace{K(K(\dots(K}}(X)\dots),
\]
where $K(X)$ denotes the degreewise free $K$-vector space on $X$. 
The $K$-completion of $X$ is derived from the cosimplicial space $R_K$ and can be defined in at least three different ways, taking values in different categories that each come with a specific notion of homotopy.

\medskip

\begin{tabular}{ccc}
Completion functor & taking values in & with homotopy defined by\\
\hline
$X\hat{{}_K} = \Tot R_K(X)$ & spaces & (simplicial) homotopy\\
$X^{\pro}_K = \{\Tot^s R_K(X)\}_{s \geq 0}$ & towers of spaces & pro-homotopy\\
$R_K$ & cosimplicial spaces & external homotopy
\end{tabular}

\medskip

For more precise information on pro-homotopical structures, see \cite{isaksen:strictmodel,fausk-isaksen:model-pro}. 
There are natural transformations $R_K \to (-)^{\pro}_K \to (-)\hat{{}_K}$ under which homotopies map to homotopies, but the three notions of homotopy are strictly contained in one another. 
To illustrate this, consider the Bousfield (unstable Adams) spectral sequence associated to $R_K$. 
If the cosimplicial space $R_K(X)$ is externally homotopy equivalent \cite{quillen:homalg,bousfield:cosimplicial-resolutions} to $R_K(Y)$, the spectral sequences agree from $E^2$ on. 
If the towers are pro-weakly equivalent, then for every $(p,q)$ there is an $R$ such that the $E^r_{p,q}$-terms are isomorphic for $r\geq R$. 
If we merely know that $X\hat{{}_K} \simeq Y\hat{{}_K}$, even if this is induced by a map $X \to Y$, the spectral sequences might not be isomorphic anywhere at any term.

The functors $(-)^{\pro}_K$ and $R_K$ can be made into endofunctors by enlarging the source category to towers or cosimplicial spaces, respectively, applying the completion functor levelwise, and taking the diagonal.

In \cite[I.5.6]{yellowmonster}, Bousfield and Kan claimed that $(-)^{\pro}_K$ is a monad (levelwise, not just pro-). In fact, they try to prove that $R_K$ is a monad. Both claims were not proved and later retracted by the authors. In the companion paper \cite{bauer-libman:monads} we showed that $(-)\hat{{}_K}$, while not a strict monad, can be made into a monad up to coherent homotopy, an \emph{$A_\infty$-monad}. This works in great generality: whenever given a monad ($A_\infty$ or strict), there is a completion functor which is again an $A_\infty$-monad.

In this paper, we address the problem of realizing this structure on the resolution. While related in its results, the methods we use here are completely different from \cite{bauer-libman:monads}, and this paper can thus be read independently.

Our result is slightly more general than outlined above.
Let $c\E$ denote the category of cosimplicial objects in a complete and cocomplete monoidal category $(\E,\comp,I)$. 
The category $c\E$ is endowed with an ``external'' simplicial structure (see Section~\ref{ringadsection}) and inherits a levelwise monoidal structure $\comp$ from $\E$.
Every object $X\in c\E$ gives rise to an operad $\opEND(X)$ (of simplicial sets) whose $n$th space is $\map_{c\E}(X^{\comp n},X)$. As usual \cite{may:72}, an \emph{$\O$-algebra} in $(c\E,\comp,I)$ is an object $X \in c\E$ with a morphism of operads $\O \to \opEND(X)$. 
For example, if $\O(n)=*$ is the associative operad, $X$ is simply a monoid.
When the spaces $\O(n)$ are weakly contracible we call $X$ an $A_\infty$-monoid.

Now let $(K,\eta,\mu)$ be a monoid in $(\E,\comp,I)$ defined by $\eta \colon I \to K$ and $\mu \colon K\comp K \to K$.
Let $R_K \in c\E$ be its cobar construction with $R_K{}^n=K^{\comp {n+1}}$ and with coface maps $d^i=K^{ \comp i}\eta K^{\comp n-i}$ and codegenracy maps $s^i=K^{\comp i}\mu K^{\comp n-i}$, cf. \cite[Ch I, \S 4]{yellowmonster}.

\begin{defn}\label{twistmap}
A \emph{twist map} for $K$ is a morphism $t \colon K \comp K \to K \comp K$ such that
\begin{enumerate}
\item \label{twmap:comm}
$\mu \circ t = \mu$,
\item \label{twmap:2}
$t \circ (K \comp \eta)= \eta \comp K$  and $t \circ (\eta \comp K) = K \comp \eta$,
\item \label{twmap:3}
$t \circ (\mu \comp K) = (K \comp \mu) \circ (t \comp K) \circ (K \comp t)$ and
$t \circ (K \comp \mu) = (\mu \comp K) \circ (K \comp t) \circ (t \comp K)$.
\end{enumerate}
\end{defn}

\begin{thm}\label{twethm}
There exists an $A_\infty$-operad of simplicial sets $\TW$ such that for every monoid $(K,\eta,\mu)$ with a twist map $t$ in a complete and cocomplete monoidal category $\E$, $R_K$ is a $\TW$-algebra in $c\E$, i.e. there is a morphism of operads $\TW \to \opEND(R_K)$.
\end{thm}

As an application, let $\E$ be the category of endofunctors of simplicial sets with composition as its monoidal structure; one needs to pass to a larger universe in order for $\E$ to make sense.
A monad $K$ on the category of spaces is a monoid in $\E$.

Let $K$ be the free $K$-module functor on spaces for a ring $K$, and $t$ the ``twist map'' defined by
\[
t = (K\comp \mu) - \id_{K \comp K} + (\mu \comp K) \quad \text{ (cf. \cite[p. 28]{yellowmonster})}
\]
We obtain as a corollary the motivating example:

\begin{cor}\label{twcomp}
For any ring $K$, $R_K$ is a $\TW$-algebra and in particular an $A_\infty$-monad.
\end{cor}

We do not know any interesting example of a monad with a twist map that is not of the free-$K$-module type. It would be interesting to know if the monad $X \mapsto \Omega^\infty(E \wedge X)$, for an $\SS$-algebra $E$ in the sense of \cite{ekmm} or \cite{hss:symmetric} can be given a twist map -- this would imply that the Bendersky-Thompson resolution \cite{bendersky-thompson} is an $A_\infty$-monad. Despite the lack of other examples, we decided to work in this generality as, if anything, it makes the proofs clearer.

%
%
%
%
%
\section{\texorpdfstring{The ``Twist map Operad'' $\TW$}{The ``Twist map Operad'' TW}} \label{simplainftyoperad}

In this section we will give an explicit combinatorial construction of the operad $\TW$ of spaces. In the next section we will construct its action on the Bousfield-Kan resolution \cite{yellowmonster}.

\subsection{Informal description} \label{informal-tw}
The $k$-simplices of the $n$th space of the operad $\TW$ are braid-like diagrams as illustrated in 
Figure~\ref{braidexample} below.
These braid diagrams have $n(k+1)$ incoming strands labelled $(0\dots k\, 0\dots k\, \dots \dots 0\dots k)$ and $k+1$ outgoing strands
labelled $(0\dots k)$ in this order from left to right.
Two strands can join if they have the same label, or they can cross provided the label of the strand on the left is
bigger than the one of the strand on the right.
In addition we allow new strands of any label to emerge.
\begin{equation}\label{crossjoin}
\text{Cross $(a>b)$:} \quad 
\begin{matrix}
\xy
(0,0)*!{a}!C; (4,-8)*!{a}!C**\dir{-}; (4,0)*!{b}!C; (0,-8)*!{b}!C**\dir{-}
\endxy
\end{matrix}
\qquad \text{Join:} \quad
\begin{matrix}
\xy
(0,0)*!{a}!C; (2,-4)**\dir{-}; (4,0)*!{a}!C**\dir{-}; (2,-4);(2,-8)*!{a}!C**\dir{-}
\endxy
\end{matrix}
\qquad \text{Emerge:} \quad
\begin{matrix}
\xy
(0,0)*+<2pt>{}*\cir{};(0,-4)*!{a}!C**\dir{-}
\endxy
\end{matrix}
\end{equation}

Note that a crossing has no orientation, that is, none of the strands is considered passing above or below the other. 
Figure~\ref{braidexample} shows an example of such a braid.
\begin{figure}[ht]
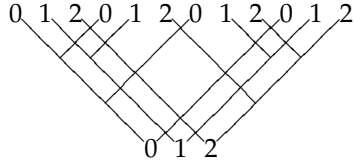

\[
\xy
(0,0)*{0};(18,-18)*{0}**\dir{-};(36,0)*{0}**\dir{-},(4,0)*{1};(22,-18)*{1}**\dir{-};(40,0)*{1}**\dir{-},(8,0)*{2};(26,-18)*{2}**\dir{-};(44,0)*{2}**\dir{-},
(12,0)*{0};(6,-6)**\dir{-},(16,0)*{1};(10,-6)**\dir{-},(20,0)*{2};(32,-12)**\dir{-},(24,0)*{0};(12,-12)**\dir{-},(28,0)*{1};(34,-6)**\dir{-},(32,0)*{2};(38,-6)**\dir{-}
\endxy
\]
\caption{An example of an element of $\TW(4)_2$} \label{braidexample}
\end{figure}
Two such diagrams are considered equal if they are isotopic in $\RR^2$ relative to the end points; in particular, we do not allow Reidemeister moves. However, we do identify the following diagrams involving the new strand operation:
\begin{equation}
\label{illustrate:equiv}
\begin{matrix}
\xy
(4,-4)*+<2pt>{}*\cir{};(8,-8)*{a}**\dir{-};(16,0)*{a}**\dir{-}
\endxy
\end{matrix} \sim
\begin{matrix}
\xy
(0,0)*{a};(0,-8)*{a}**\dir{-}
\endxy
\end{matrix} \sim
\begin{matrix}
\xy
(12,-4)*+<2pt>{}*\cir{};(8,-8)*{a}**\dir{-};(0,0)*{a}**\dir{-}
\endxy
\end{matrix}\qquad,\quad
\begin{matrix}
\xy
(2,-2)*+<2pt>{}*\cir{};(8,-8)*{a}**\dir{-};(0,-8)*{b};(8,0)*{b}**\dir{-}
\endxy
\end{matrix}\sim
\begin{matrix}
\xy
(0,0)*{b};(0,-8)*{b}**\dir{-};(4,-4)*+<2pt>{}*\cir{};(4,-8)*{a}**\dir{-}
\endxy
\end{matrix}\qquad,\quad
\begin{matrix}
\xy
(6,-2)*+<2pt>{}*\cir{};(0,-8)*{b}**\dir{-};(8,-8)*{a};(0,0)*{a}**\dir{-}
\endxy
\end{matrix}\sim
\begin{matrix}
\xy
(4,0)*{a};(4,-8)*{a}**\dir{-};(0,-4)*+<2pt>{}*\cir{};(0,-8)*{b}**\dir{-}
\endxy
\end{matrix}
\end{equation}
The sequence of $\{\TW(n)_k\}_{k \geq 0}$ forms a simplicial set, where the $i$th face map $d_i$ is obtained by erasing the strands of color $i$, whereas the $i$th degeneracy map $s_i$ duplicates the $i$th strand. 
An illustration of how this works on the ``basic building blocks'' of the braids is depicted by (here $a>b$)
\begin{equation}
\label{illustrate:dup}
s_a
\left(\begin{matrix}
\xy
(0,0)*!{a}!C; (3,-6)**\dir{-}; (6,0)*!{a}!C**\dir{-}; (3,-6);(3,-12)*!{a}!C**\dir{-}
\endxy
\end{matrix}\right)
=\begin{matrix}\xy
(0,0)*!{a}!C; (6,-6)**\dir{-}; (12,0)*!{a}!C**\dir{-}; (6,0)*!{a+1}!C; (12,-6)**\dir{-}; (18,0)*!{a+1}!C**\dir{-};(6,-6);(6,-12)*!{a}!C**\dir{-};(12,-6);(12,-12)*!{a+1}!C**\dir{-}
\endxy\end{matrix}
,\qquad
s_b
\left(\begin{matrix}
\xy
(0,0)*!{a};(6,-12)*!{a}**\dir{-},
(6,0)*!{b};(0,-12)*!{b}**\dir{-}
\endxy
\end{matrix}\right)
=\begin{matrix}\xy
(-2,0)*!{a+1}!C;(16,-12)*!{a+1}!C**\dir{-},
(6,0)*!{b}!C;(0,-12)*!{b}!C**\dir{-},
(12,0)*!{b+1}!C;(6,-12)*!{b+1}!C**\dir{-}
\endxy\end{matrix}
\end{equation}
thus introducing new crossings.

Composition in the operad $\TW$ is obtained by grafting of these braid-like diagrams with matching labels.

We will now give a more rigid combinatorial description of the operad $\TW$ without direct reference to the braid-like diagrams.
To do this we will construct a simplicial object $\{\B_k\}_k$ in the category of small monoidal categories. The $k$-simplices of the space $\TW(n)$ are the morphisms between certain objects in $\B_k$.
In fact, each $\B_k$ is a free monoidal category in a sense which we will now make precise.

%
\subsection{\texorpdfstring{Free monoidal categories and the simplicial monoidal category $\{\B_k\}_{k \geq 0}$}{Free monoidal categories and the simplicial monoidal category B}}

By a \emph{graph} we mean a pair $(\A \rightrightarrows \O)$ of two sets and two maps (codomain, domain) between them. A \emph{unital graph} is a graph with a common section $\O \to \A$ of domain and codomain. 
Every small category $\C$ has an underlying unital graph $U_u(\C) = (\Mor\C \overset{\curvearrowleft}{\rightrightarrows} \Obj \C)$ by forgetting the composition rule on morphisms.
By further forgetting the unit maps, we obtain a graph $U(\C)$. Conversely, every (unital) graph $\G = ( \A \rightrightarrows \O)$ has an associated category $C(\G)$. It is the free category generated by $\G$ in the sense that $C$ is left adjoint to the functor $U$ (resp. $U_u$) from the category of small categories to the category of (unital) graphs. See \cite[Section II.7, in particular Theorem~1]{maclane:cwm}. 
This section is about a version of this free construction for monoidal categories.

\begin{defn}
A \emph{vertex-monoidal graph} is a graph $(\A \rightrightarrows \O)$, where $\O$ has the structure of a monoid. A \emph{(unital) monoidal graph} is a graph object in the category of (unital) monoids. 
\end{defn}

\begin{lem}
\label{free-monoidal}
The following forgetful functors all have left adjoints:
\begin{align*}
\{\text{small monoidal categories}\} \xto{U_{um}} & \{\text{unital monoidal graphs}\}\\
\xto{U_m} & \{\text{monoidal graphs}\} \xto{U_{vm}} \{\text{vertex-monoidal graphs.}\}
\end{align*}
In particular, for every vertex-monoidal graph $\G$, there exists a small monoidal category $(\M,\sqcup,\emptyset)$ with a morphism of vertex-monoidal graphs $\epsilon\colon \G \to U(\M)$, such that for any other monoidal category $(\C,\otimes,I)$ and vertex-monoidal graph morphism $X\colon \G \to U(\C)$, there is a unique functor $X^{\flat} \colon \M \to \C$ of monoidal categories such that $\epsilon \circ U(X^\flat) = X$.
\end{lem}
We give a proof of this (easy) fact in the appendix.

\begin{remark} \label{rem-free-monoidal}
Fix a set $\Sigma$, and let $\Sigma_{\sqcup}$ be the set of words, i.e. finite sequences in $\Sigma$, including the empty sequence $()$. 
This is a monoid under concatenation, which we denote by $\sqcup$. 
If $(\A \rightrightarrows \O)$ is a vertex-monoidal graph whose vertices $\O = \Sigma_\sqcup$ form such a free monoid, then a morphism of vertex-monoidal graphs $X\colon (\A \rightrightarrows \Sigma_\sqcup) \to U(\M)$, where $\M$ is monoidal, is given by a map $X\colon \Sigma \to \Obj(\M)$ and, for every arrow $a\colon (\sigma_1,\dots,\sigma_k) \to (\tau_1,\dots,\tau_l)$ in $\A$, a morphism $X(\sigma_1) \otimes \dots \otimes X(\sigma_k) \to X(\tau_1) \otimes \dots \otimes X(\tau_l)$.
\end{remark}

\begin{lem}
\label{free-monidel-w-relations}
Let $\Sigma$ and $\G=(\A \rightrightarrows \Sigma_\sqcup)$ be as in Remark~\ref{rem-free-monoidal} and let $\M$ denote the free monoidal category generated by $\G$ by Lemma~\ref{free-monoidal}.
Let $\sim$ be an equivalence relation on the morphism set of $\M$ such that if $f \sim f'$ then their domains and codomains are equal.

Then there exists a monoidal category $\M/_\sim$ together with a quotient functor $\M \to \M/_\sim$ which is the identity on the object sets and such that $\M/_\sim$ is the free monoidal category on $(\Sigma,\Lambda)$ subject to the relation $\sim$.
That is, given any  monoidal category $(\C,\otimes,I)$ and a morphism of monoidal graphs 
$X^\# \colon \G \to U(\C)$ whose associated monoidal functor $X \colon \M \to \C$ has the property that $X(f)=X(g)$ if $f \sim g$, there exists a unique monoidal functor $X \colon \M/_\sim \to \C$ which factors $X$ through the quotient $\M \to \M/_\sim$.
\end{lem}

\begin{proof}
Consider the equivalence relation $\approx$ on the morphism set of $\M$ which is generated by $\sim$  and which is closed to the property that $(f \sqcup g) \approx (f' \sqcup g')$ and  $f \circ g \approx f' \circ g'$ if $f \approx f'$ and $g \approx g'$ (for the composition, whenever it makes sense).
One easily checks that the domain and codomain of $\approx$-equivalent morphisms are equal.
We can now apply \cite[\S II.8, Proposition 1]{maclane:cwm} to obtain $\M/_\approx$ which is monoidal due to the requirements on $\approx$.
\end{proof}

The definition of the category $\B_k$ is modelled on the ideas in \ref{informal-tw}.

\begin{defn}\label{def-Bk}
We use Lemma \ref{free-monidel-w-relations} to define $\B_k$ as the free monoidal category generated by the set of symbols
 $[k]=\{0,\ldots,k\}$ and morphisms (cf.~\eqref{crossjoin})
\begin{eqnarray*}
() \xto{u^a} (a) \qquad \text{ and } \qquad 
(a,a) \xto{m^a} (a) & & \text{ for all } a=0,\ldots,k \\
(a,b) \xto{t^{a,b}} (b,a) & & \text{ for all } 0 \leq b<a  \leq k
\end{eqnarray*}
subject to the relations (cf.~\ref{illustrate:equiv})
\begin{eqnarray*}
(R1)           & & m^a \circ (\id_{(a)} \sqcup u^a) = \id_{(a)} = m^a \circ (u^a \sqcup \id_{(a)}) \\
(R2)_i         & & t^{c,d} \circ (\id_{(c)} \sqcup u^d) = u^d \sqcup \id_{(c)} \\
(R2)_{ii}     & & t^{c,d} \circ (u^c \sqcup \id_{(d)}) = \id_{(d)} \sqcup u^c 
\end{eqnarray*}
where $0 \leq a,c,d \leq k$ and $d<c$.
\end{defn}

For an object $\si=(a_1,\cdots,a_n)$  and a morphism $f$ in $\B_k$ we will write $\si \sqcup f$ for $\id_{\si} \sqcup f$.
The morphism in $\B_k$ between sequences $(a_1,\ldots,a_s) \to (b_1,\ldots,b_t)$ 
can be depicted by braid-like diagrams with $s$ incoming strands labelled by the integers $a_1,\ldots,a_s$ 
and $r$ outgoing strands labelled by $(b_1,\ldots,b_t)$.
For example, the diagram
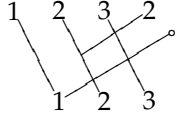
\[
\xy
(0,0)*{1};(6,-12)*{1}**\dir{-};(21,-3)*+<2pt>{}*\cir{}**\dir{-},
(6,0)*{2};(12,-12)*{2}**\dir{-},
(12,0)*{3};(18,-12)*{3}**\dir{-},
(18,0)*{2};(9,-6)**\dir{-}
\endxy
\]
describes a morphism $(1232) \to (123)$ in $\B_3$ which is the composition
\[
[m^1 \sqcup (23)] \circ [(1) \sqcup t^{21} \sqcup (3)] \circ
[(1) \sqcup m^2 \sqcup t^{31}] \circ [(12) \sqcup t^{32}\sqcup u^1].
\]

Most of this section is devoted to the proof of the following result.

\begin{prop}
\label{bk-simplcl}
There are face and degeneracy operators $d_i \colon \B_k \to \B_{k-1}$ and $s_i \colon \B_k \to \B_{k+1}$ which make
$\{ \B_k\}_{k \geq 0}$ a simplicial object of small monoidal categories.
On objects, $d_i$ and $s_i$ act by deleting (resp. duplicating) the symbol $i \in [k]$.
\end{prop}

A more elaborate description of the simplicial operators is given in Definition~\ref{Bk-simplicial}.
In what follows we will write $[k]$ for the object $(01\dots k)$ in $\B_k$ use $[k]^n$ to denote the $n$-fold concatenation of $[k]$ with itself, i.e. the object $(0\dots k 0\dots k \dots)$.

\begin{defn}\label{def-tw}
Define $\TW(n)_k = \Hom_{\B_k}([k]^n,[k])$.
By Proposition~\ref{bk-simplcl}, this collection of sets, for all $k$, forms a simplicial set.

The composition law in the operad 
$
\TW(n) \times \prod_{i=1}^n \TW(m_i) \longrightarrow \TW(\sum_{i=1}^n m_i)
$
is defined using the monoidal structure of $\B_k$:
\begin{multline*}
\B_k([k]^n,[k]) \times \prod_{i=1}^n \B_k([k]^{m_i},[k]) \xto{\quad \id \times \sqcup \quad}
\B_k([k]^n,[k]) \times \B_k([k]^{\sum_{i} m_i},[k]^n) 
\\
\xto{\quad \circ \quad} \B_k([k]^{\sum_i m_i},[k]).
\end{multline*}
Unitality and associativity of the monoidal structure in $\B_k$ descend to this composition law.
\end{defn}

The goal of this section is to prove

\begin{thm}\label{tw-is-ainfty}
The sequence $\TW(n)$ is a (nonsymmetric) $A_\infty$-operad of simplicial sets.
\end{thm}

Here are some easy results on the categories $\B_k$ (Definition \ref{def-Bk}.)
An object $\si \in \B_k$ is a sequence of symbols $(s_1,\ldots,s_n)$ and we write $u^{\si}$ for 
$u^{s_1} \sqcup \dots \sqcup u^{s_n}$.

\begin{lem}\label{initial}
Fix an object $\tau = (t_1 < \dots < t_n)$ in $\B_k$.
Then given an object $\sigma \in \B_k$, any morphism $\vp \colon \tau \to \sigma$ has the form 
\begin{equation}\label{phidecomp}
\vp = u^{\si_0} \sqcup (t_1) \sqcup u^{\si_1} \sqcup \dots \sqcup u^{\si_{n-1}} \sqcup  (t_n) \sqcup u^{\si_n}
\end{equation}
where $\si_i$ are subsequences of $\si$ such that $\si=\si_0 \sqcup (t_1) \sqcup \si_1 \sqcup \dots \sqcup \si_{n-1} \sqcup (t_n) \sqcup \si_n$. In particular, the object $()$ is initial in $\B_k$ and $\End_{\B_k}(\tau)=\{\id\}$.
\end{lem}

\begin{proof}
By Lemma \ref{free-monidel-w-relations}, the morphisms in $\B_k$ are compositions of morphisms of the form $x \sqcup u^a \sqcup y$, $x \sqcup m^a \sqcup y$, and $x \sqcup t^{ab} \sqcup y$, where $x$, $y$ are objects of $\B_k$.

Let $\vp$ be as in the statement of the lemma. We prove the result by induction on the ``generation length'' $s$ of $\vp$. If $s=0$ then $\vp=\id$.
If $\vp$ satisfies \eqref{phidecomp} and is followed by some $x \sqcup u^a \sqcup y$, then the resulting morphism still has the form \eqref{phidecomp}.
If $\vp$ is followed by some $x \sqcup t^{ab} \sqcup y$  then $(a,b)$ cannot be contained in $\tau$ because $a>b$, hence one of $a$ or $b$ belongs to $\si_i$ and we can apply the relation (R2) to show that the composition 
is again of the form \eqref{phidecomp}.
Finally, if $\vp$ is followed by $x \sqcup m^a \sqcup y$ then the domain $(a,a)$ of $m^a$ cannot be contained in $\tau$, and relation (R1) guarantees that the composition is of the form \eqref{phidecomp}.
\end{proof}

\begin{defn}
\label{Bk-simplicial}
We turn the categories $\B_k$ into a simplicial monoidal category as follows.
To do this it is convenient to replace $[k]$ with a totally ordered finite set $S$ of cardinality $k+1$.
Then monotonic functions $S \to T$ are generated by the following functions
\begin{eqnarray*}
& & \eta_a \colon (S -\{ a\}) \xto{\text{  inclusion }} S \qquad \text{ for some } a \in S \\
& & \epsilon_a \colon (S-\{a\}) \sqcup \{ a',a''\} \xto{\text{  collapse } \{a',a''\} \to \{a\}} S 
\end{eqnarray*}
where $a'<a''$ are inserted to $S$ instead of $a \in S$.
In this notation the simplicial identities become
\begin{equation}\label{simpidntsee:1}
\begin{matrix}
\eta_a \circ \eta_b = \eta_b \circ \eta_a, &
    \epsilon_b \circ \eta_a = \eta_a \circ \epsilon_b, &
    \epsilon_a \circ \epsilon_b = \epsilon_b \circ \epsilon_a  & (a \neq b \in S) \\
& \epsilon_a \circ \epsilon_{a'} = \epsilon_a \circ \epsilon_{a''}, 
    & \epsilon_a \circ \eta_{a'} = \epsilon_a \circ \eta_{a''} =\Id  &  (a \in S) 
\end{matrix}
\end{equation}
We define face operators $d_a = (\eta_a)^* \colon \B_S \to \B_T$ where $T=S-\{ a \}$ and  degeneracy operators
$s_a = (\epsilon_a)^*\colon \B_T \to \B_S$ where $T=S-\{a\} \sqcup\{ a'<a''\}$ as follows.
On objects we define 
\[
d_a((b))=(b) = s_a((b)) \text{ if $b \neq a$;} \quad d_a((a))=(); \quad  s_a((a))=(a',a'')
\]
Extend $s_a$ and $d_a$ to all the objects of $\B_k$ by monoidality with respect to $\sqcup$.

On morphisms we define $d_a$ and $s_a$ by defining them on generators:
$$
\begin{array}{llr}
d_a(m^a)=\id_{()}, & s_a(m^a)= (m^{a'}\sqcup m^{a''}) \circ ((a') \sqcup t^{a'',a'} \sqcup (a''))& \text{cf. \eqref{illustrate:dup}}\\
d_a(m^b)=m^b, & s_a(m^b)=m^b & \text{if } b \neq a \\
d_a(u^a)=\id_{()}, & s_a(u^a)=u^{a'} \sqcup u^{a''}\\
d_a(u^b)=u^b, & s_a(u^b)=u^b & \text{if } b \neq a \\
d_a(t^{cd})= t^{cd}, & s_a(t^{cd}) = t^{cd} & \text{if } c,d \neq a \\
d_a(t^{cd})= \id_{(d)}, &  s_a(t^{cd}) =(t^{a' d} \sqcup (a'')) \circ ((a') \sqcup t^{a'' d}) & \text{if } c=a\\
d_a(t^{cd})= \id_{(c)}, &  s_a(t^{cd}) =	((a') \sqcup t^{c a''}) \circ (t^{ca'} \sqcup (a'')) & \text{if } d=a, \text{cf. \eqref{illustrate:dup}}\\
\end{array}
$$
To see that these assignments define functors we need to check that they respect the relations (R1) and (R2).
For (R1), we compute ($a \neq b$):
\begin{eqnarray*}
s_b(m^a \circ ((a) \sqcup u^a)) &=& m^a \circ ((a) \sqcup u^a) = \id_{(a)} = s_b(\id_{(a)}) \\
s_a(m^a \circ ((a) \sqcup u^a)) &=& 
  (s^{a'} \sqcup s^{a''}) \circ ((a') \sqcup t^{a''a'} \sqcup (a'')) \circ ((a'a'') \sqcup u^{a'}\sqcup u^{a''}) \\
 & \underset{\text{(R2)}} =  & (s^{a'} \sqcup s^{a''}) \circ ((a') \sqcup u^{a'} \sqcup (a'') \sqcup u^{a''})\\
 &  \underset{\text{(R1)}} = & \id_{(a')} \sqcup \id_{(a'')} =
\id_{(a'a'')} =s_a(\id_{(a)}).
\end{eqnarray*}
Similarly one checks that $d_a$ and $d_b$ carry the relation (R1) to equal morphisms.
The second relation of (R1) is checked similarly.

For (R2), assume that $c>d$ and consider some $a$.
If $a \neq c,d$ then it is clear that applying $d_a$ and $s_a$ to the relations (R2) leaves them unaffected.
If $a=c$ then applying $s_a$ to the relation (R2)${}_{i}$ gives
\begin{multline*}
s_a (t^{ad} \circ ((a) \sqcup u^d)) = 
    (t^{a'd} \sqcup (a'')) \circ ((a') \sqcup t^{a'd}) \circ ((a'a'') \sqcup u^d) \underset{(R2)}{=} \\
(t^{a'd} \sqcup (a'')) \circ ((a') \sqcup u^d \sqcup(a'')) \underset{(R2)}{=} 
u^d \sqcup \id_{(a'a'')} = 
    s_a(u^d \sqcup \id_{(a)}) 
\end{multline*}
When $a=d$ a similar calculation shows that $s_a(t^{ca} \circ (\id_{(c)} \sqcup u^a)) =  s_a(u^a \sqcup \id_{(c)})$. 
Checking that $s_a$ respects the relation $(R2)_{ii}$ is analogous, and we leave it to the reader to check that $d_a$ respects (R2) as well. \qed
\end{defn}

\begin{proof}[Proof of Proposition \ref{bk-simplcl}]
We will use the ``coordinate free'' notation of Definition~\ref{Bk-simplicial} and prove that the functors $d_a$ and $s_a$ render $\{\B_k\}_{k \geq 0}$ a simplicial monoidal category.
The functors $d_a$ and $s_a$ are monoidal, so it only remains to prove that they satisfy the simplicial identities dual to the identities \eqref{simpidntsee:1}, i.e.
\begin{eqnarray}
d_b d_a &=& d_a d_b \quad (a \neq b \in S) \label{dbda}\\
d_b s_a &=& s_a d_b \quad (a \neq b \in S)\label{dbsa}\\
s_b s_a &=& s_a s_b \quad (a \neq b \in S)\label{sbsa}\\
s_{a'} s_a &=& s_{a''} s_a \quad (a \in S)\label{sasa}\\
d_{a'}s_a &=& d_{a''} s_a = \id \quad (a \in S)\label{dasa}
\end{eqnarray}
This is obviously the case on the object sets and it only remains to check the simplicial identities on the generating 
morphisms $u^c, m^c, t^{cd}$.
If $c,d \in S - \{ a,b\}$ then by definition $d_a, d_b$ as well as $s_a$ and $s_b$ map $u^c, m^c$ and $t^{cd}$ to
themselves, hence all the identities hold.
Thus, it remains to check the simplicial identities when $\{c,d\}$ and $\{a,b\}$ are not disjoint.
For \eqref{dbda} and \eqref{sbsa}, the only non-immediate cases are
\begin{eqnarray*}
s_b s_a (m^a) =& s_b (m^{a'} \sqcup m^{a''}) \circ s_b ({(a')} \sqcup t^{a''a'} \sqcup {(a'')})\\
 =& (m^{a'} \sqcup m^{a''}) \circ ({(a')} \sqcup t^{a''a'} \sqcup {(a'')}) &= s_a s_b (m^a)\\
s_b s_a (t^{ab}) =& s_b (t^{a'b} \sqcup {(a'')}) \circ s_b({(a')} \sqcup t^{a''b})\\
=& ({(b')} \sqcup t^{a'b''} \sqcup {(a'')}) \circ (t^{a'b'} \sqcup {(b'')} \sqcup {(a'')})\\
 & \circ({(a')} \sqcup {(b')} \sqcup t^{a''b''}) \circ ({(a')} \sqcup t^{a''b'} \sqcup {(b'')})\\
=& ({(b')} \sqcup t^{a'b''} \sqcup {(a'')}) \circ (t^{a'b'} \sqcup t^{a''b''}) \circ ({(a')} 
    \sqcup t^{a''b'} \sqcup {(b'')})\\
=& ({(b')} \sqcup t^{a'b''} \sqcup {(a'')}) \circ ({(b')} \sqcup {(a')} \sqcup t^{a''b''})\\
 & \circ (t^{a'b'} \sqcup {(a'')} \sqcup {(b'')}) \circ ({(a')} \sqcup t^{a''b'} \sqcup {(b'')})\\
=& s_a ({(b')} \sqcup t^{ab''}) \circ s_a (t^{ab'} \sqcup {(b'')}) &= s_a s_b (t^{ab})\\
\end{eqnarray*}
The remaining cases, for example $d_b d_a (m^a) = d_a d_b (m^a)$ or when $\{a,b\} \neq \{c,d\}$, are straightforward.

For \eqref{dbsa}, we check
\begin{eqnarray*}
d_b s_a (m^a) =& d_b (m^{a'} \sqcup m^{a''}) \circ d_b({(a')} \sqcup t^{a''a'} \sqcup {(a'')})\\
 =& (m^{a'} \sqcup m^{a''}) \circ ({(a')} \sqcup t^{a''a'} \sqcup {(a'')}) &= s_a d_b (m^a)\\
d_b s_a (t^{ac}) =& d_b (t^{a'c} \sqcup {(a'')}) \circ d_b ({(a')} \sqcup t^{a''c})\\
=& (t^{a'c} \sqcup {(a'')}) \circ ({(a')} \sqcup t^{a''c}) &= s_a d_b (t^{ac})\\ 
\end{eqnarray*}
and leave the rest to the reader.

For \eqref{sasa}, we use the convention that if $a \in S$ then the domain of $\epsilon_a \colon S' \to S$ is formed from $S$ by replacing $a$ with $a'< a''$. 
The morphisms $\epsilon_{a'},\epsilon_{a''} \colon S'' \to S'$ have $S''=S-\{a\} \sqcup \{a',a'',a'''\}$
and $s_{a'}(a') = (a',a'')$, $s_{a'}(a'')=a'''$ whereas  $s_{a''}(a')=a'$, and $s_{a''}(a'')=(a'',a''')$. We now have the following calculation, which exploits the monoidal structure of $\B_k$:
\begin{eqnarray*}
s_{a'} s_{a} (m^a) =& 
  s_{a'}(m^{a'} \sqcup m^{a''}) \circ s_{a'}((a') \sqcup t^{a''a'} \sqcup (a''))  \\
  =& \{\lbrack (m^{a'} \sqcup m^{a''}) \circ((a') \sqcup t^{a''a'} \sqcup (a'')) \rbrack \sqcup m^{a'''}\}\\
   &  \circ\{ (a'a'') \sqcup \lbrack ((a') \sqcup t^{a'''a''}) \circ (t^{a'''a'} \sqcup (a'')) \rbrack \sqcup (a''') \} \\
  =& (m^{a'} \sqcup m^{a''} \sqcup m^{a'''}) \circ ((a') \sqcup t^{a''a'} \sqcup (a''a'''a''')) \\
   &  \circ((a'a''a')\sqcup t^{a'''a''}\sqcup (a''')) \circ ((a'a'') \sqcup t^{a'''a'} \sqcup (a''a''')) \\
  =& (m^{a'} \sqcup m^{a''} \sqcup m^{a'''}) \circ ((a'a'a'') \sqcup t^{a'''a''} \sqcup (a''')) \\
   & \circ ((a')\sqcup t^{a''a'}\sqcup (a'''a''a''')) \circ ((a'a'') \sqcup t^{a'''a'} \sqcup (a''a''')) \\
  =& \{ m^{a'} \sqcup \lbrack (m^{a''} \sqcup m^{a'''}) \circ ((a'') \sqcup t^{a'''a''} \sqcup (a''')) \rbrack \}\\
   &  \circ \{ (a') \sqcup \lbrack (t^{a''a'} \sqcup (a''')) \circ ((a'') \sqcup t^{a'''a'}) \rbrack \sqcup (a''a''') \} \\
  =& s_{a''}(m^{a'} \sqcup m^{a''}) \circ s_{a''}((a') \sqcup t^{a''a'} \sqcup (a'')) &= s_{a''} s_{a}(m^a) 
\end{eqnarray*}
This is the most subtle calculation in verifying \eqref{dbsa}.
The remaining cases, e.g. $s_{a'}s_a(u^a) = s_{a''}s_a(u^a)$, are proved similarly and left to the reader, as is the proof of \eqref{dasa}.
\end{proof}

\begin{prop} \label{twzeroonepoint}
The spaces $\TW(0)$ and $\TW(1)$ are points.
\end{prop}
\begin{proof}
This follows immediately from Lemma~\ref{initial}.
\end{proof}

\begin{remark}
\label{cofacial}
The categories $\B_k$ also enjoy the structure of a ``cofacial monoidal category'', namely a cosimplicial object in the
category of small monoidal categories without codegeneracy maps.
In the notation of Definition~\ref{Bk-simplicial} it is induced by the inclusion of symbols $T \subseteq S$ where 
$T$ and $S$ are totally ordered sets.
We will not need this structure but we will make use of the operator 
\[
d^0 \colon \B_k \to \B_{k+1}
\]
which has the effect of shifting indices by $1$.
That is, $d^0(i)=i+1$ for all $i \in [k]$ and on morphisms $d^0(u^a)=u^{a+1}, d^0(m^a)=m^{a+1}$ and 
$d^0(t^{ab})=t^{a+1,b+1}$.
It is obvious that the relations (R1) and (R2) are respected by these assignments.
The ``coface operator'' $d^0$ also has the following relations with the face and degeneracy operators:
\begin{equation}\label{cofacialid}
\begin{array}{ccc}
d_0 d^0 = \id, & s_0 d^0 = d^0 d^0 \\
d_i d^0 = d^0 d_{i-1}, & s_i d^0 = d^0 s_{i-1} & \text{ for all } i>0.
\end{array}
\end{equation}
\end{remark}

If $\sigma$ is some sequence, we use the power notation $\sigma^n$ for the $n$-fold concatenation of $\sigma$ with itself. We observe that $d_a([k+1]^n)=[k]^n=s_a([k-1]^n)$ and $d^0[k]=(1,2,\ldots,k+1)$.

\begin{defn}\label{be-n1}
For every $n \geq 2$ let $b_n \in \B_1$ denote the object $(0 1^n)$.
Define $\be_{n} \colon [1]^n \to b_n$ by induction on $n\geq 2$ as follows:
\[
\begin{matrix}
\be_2 =\{m^0 \sqcup(11)\} \circ \{ (0) \sqcup t^{10} \sqcup (1)\}, &
\be_n = \{ \be_{n-1} \sqcup (1) \} \circ \{ [1]^{n-2} \sqcup \be_2\} \\
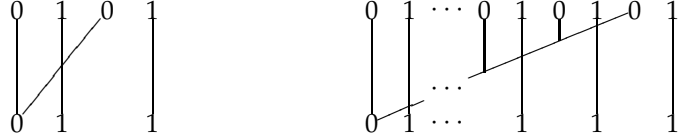
\xy
(0,0)*{0};(0,-15)*{0}**\dir{-};(12,0)*{0}**\dir{-},
(6,0)*{1};(6,-15)*{1}**\dir{-},
(18,0)*{1};(18,-15)*{1}**\dir{-}
\endxy
&
\xy
(0,0)*{0};(0,-15)*{0}**\dir{-};(7,-12)**\dir{-},
(13,-9);(35,0)*{0}**\dir{-},
(10,-10.5)*{\cdots},
(10,0)*{\cdots},
(10,-15)*{\cdots},
(5,0)*{1};(5,-15)*{1}**\dir{-},
(15,0)*{0};(15,-8.3)**\dir{-},
(20,0)*{1};(20,-15)*{1}**\dir{-},
(25,0)*{0};(25,-4.1)**\dir{-},
(30,0)*{1};(30,-15)*{1}**\dir{-},
(40,0)*{1};(40,-15)*{1}**\dir{-}
\endxy
\end{matrix}
\]
\end{defn}

\begin{lem}\label{bengames}
The following hold for all $n \geq 2$.
\begin{enumerate}
\item
\label{bengames:1}
$d_0(\be_n) = \id_{[0]^n} \in \B_1$.

\item
\label{bengames:3}
$s_0(\be_n) = \{ (0) \sqcup d^0(\be_n)\} \circ s_1 \be_n$ \quad (in $\B_2$).

\item
\label{bengames:4}
$s_0 d_1 (\be_n) =\{ (0) \sqcup d^0d_1(\be_n)\} \circ \be_n$ \quad (in $\B_2$).
\end{enumerate}
\end{lem}

The content of \eqref{bengames:3} is displayed as the obvious equivalence of the following ``braids''.
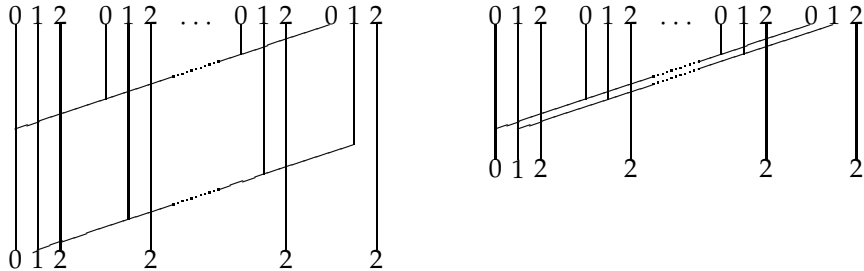
\[
\xy
(0,0)*!{0};(0,-30)*!U{0}!U**\dir{-},
(3,0)*!{1};(3,-30)*!U{1}!U**\dir{-},
(6,0)*!{2};(6,-30)*!U{2}!U**\dir{-},
(12,0)*!{0};(12,-10)**\dir{-},
(15,0)*!{1};(15,-26)**\dir{-},
(18,0)*!{2};(18,-30)*!U{2}!U**\dir{-},
(24,0)*{\cdots},
(30,0)*!{0};(30,-4)**\dir{-},
(33,0)*!{1};(33,-20)**\dir{-},
(36,0)*!{2};(36,-30)*!U{2}!U**\dir{-},
(42,0)*!{0};(27,-5)**\dir{-};
(21,-7)**\dir{.};
(0,-14)**\dir{-},
(45,0)*!{1};(45,-16)**\dir{-};(27,-22)**\dir{-};
(21,-24)**\dir{.};
(3,-30)**\dir{-},
(48,0)*!{2};(48,-30)*!U{2}!U**\dir{-}
\endxy
\qquad  \qquad
\xy
(0,0)*!{0};(0,-18)*!U{0}!U**\dir{-},
(3,0)*!{1};(3,-18)*!U{1}!U**\dir{-},
(6,0)*!{2};(6,-18)*!U{2}!U**\dir{-},
(18,0)*!{2};(18,-18)*!!U{2}!U**\dir{-},
(24,0)*{\cdots},
(36,0)*!{2};(36,-18)*!U{2}!U**\dir{-},
(48,0)*!{2};(48,-18)*!U{2}!U**\dir{-},
(12,0)*!{0};(12,-10)**\dir{-},
(15,0)*!{1};(15,-10)**\dir{-},
(30,0)*!{0};(30,-4)**\dir{-},
(33,0)*!{1};(33,-4)**\dir{-},
(42,0)*!{0};(27,-5)**\dir{-};(21,-7)**\dir{.};(0,-14)**\dir{-},
(45,0)*!{1};(27,-6)**\dir{-};(21,-8)**\dir{.};(3,-14)**\dir{-}
\endxy
\]
The content of \eqref{bengames:4} is obtained by deleting all the strands labelled with $2$.

\begin{proof}
\eqref{bengames:1}
By inspection of Definitions~\ref{Bk-simplicial} and \ref{be-n1},
$$
d_0(\be_2) = d_0(m^0\sqcup (11)) \circ d_0((0)\sqcup t^{10}\sqcup (1)) = \id.
$$
The result follows easily using the inductive definition of $\be_n$.

\eqref{bengames:3}
Using Definition~\ref{be-n1} and the pictures given there, we first observe that 
$\be_n = \ga_n \sqcup (1) $ for some  $\ga_n \colon [1]^{n-1} \to b_{n-1}$. It now follows, using the monoidal structure in $\B_2$, that for every morphism $\vp\colon [1]^2 \to b_2$ in $\B_1$ we have
\begin{multline}
\label{bengames:p1}
\{ s_1b_{n-1} \sqcup d^0\vp\} \circ \{s_1\be_n \sqcup(12)\} = 
 \{ s_1b_{n-1} \sqcup d^0\vp\} \circ \{s_1\ga_n \sqcup(12)^2\} \\
=\{ s_1\ga_n \sqcup (122) \} \circ \{ [2]^{n-1} \sqcup (0) \sqcup d^0 \vp\} =
  \{s_1\be_n \sqcup (2)\} \circ \{ [2]^{n-1} \sqcup (0) \sqcup d^0\vp\}.
\end{multline}
To complete the proof of \eqref{bengames:3} we use induction on $n$.
The hardest part is the base of induction $n=2$.
Using Definition~\ref{be-n1} and the monoidal structure in $\B_2$ we find
\begin{multline*}
\{ (0) \sqcup d^0\be_2\} \circ s_1(\be_2) = \\
\underbrace{ 
    \{ (0) \sqcup m^1 \sqcup(22)\} \circ  \{(01)\sqcup t^{21} \sqcup(2) \} \circ \{ m^0 \sqcup(1212)\}}_{(0)\circ m^0 = m^0\circ (00)} 
\circ
    \underbrace{ \{ (0) \sqcup s_1(t^{10}) \sqcup (12)\}}_{\text{expand}} = \\
\{ m^0 \sqcup m^1 \sqcup (22)\} \circ 
    \underbrace{ \{(001) \sqcup t^{21} \sqcup (2)\} \circ \{(0)\sqcup t^{10} \sqcup (212)\}} \circ
    \{(01) \sqcup t^{20} \sqcup (12)\} = \\
\{ m^0 \sqcup m^1 \sqcup (22)\} \circ \{(0) \sqcup t^{10} \sqcup (122)\} \circ
    \{(010) \sqcup t^{21} \sqcup (2)\} \circ \{(01) \sqcup t^{20} \sqcup (12)\}  = \\
s_0(m^0 \sqcup (11)) \circ s_0((0) \sqcup t^{10} \sqcup (1)) = s_0(\be_2).
\end{multline*}
Now assume that the formula holds for $n$ and we prove it for $\be_{n+1}$ by using its inductive definition.
Note that $(0) \sqcup d^0[1]^n = (01212\dots 12)=s_1(b_n)$.
\begin{multline*}
\{ (0) \sqcup d^0\be_{n+1}\} \circ s_1\be_n) \\
= \{(0) \sqcup d^0 \be_{n} \sqcup (2)\} \circ 
    \underbrace{\{ s_1b_{n-1} \sqcup d^0\be_2\} \circ 
        \{ s_1\be_{n} \sqcup (12)\} }_{\text{Formula \eqref{bengames:p1}}}
    \circ \{[2]^{n-1} \sqcup s_1\be_2\} \\
= \underbrace{\{ (0) \sqcup d^0\be_{n} \sqcup (2)\} \circ \{ s_1\be_{n} \sqcup (2)\}}_{\text{induction}} \circ
   \underbrace{ \{[2]^{n-1} \sqcup(0)\sqcup d^0\be_2\} \circ \{ [2]^{n-1} \sqcup s_1\be_2\}}_{\text{induction}} \\
=s_0(\be_{n} \sqcup (1)) \circ \{ [2]^{n-1}\sqcup s_0\be_2\} =
s_0(\be_{n+1}).
\end{multline*}
This completes the proof of \eqref{bengames:3}. Part \eqref{bengames:4} follows by applying $d_2$ to \eqref{bengames:3}, using the simplicial identities and \eqref{cofacialid}.
\end{proof}

\begin{remark}
In Proposition \ref{twzeroonepoint} we saw that $\TW(0)=\TW(1)=*$.
One can also check that $\TW(2)$ is a point, whereas all $\TW(n)$ are infinite-dimensional for $n \geq 3$.
Since we will not use these facts, we omit the proofs.
\end{remark}


\begin{proof}[Proof of Theorem \ref{tw-is-ainfty}]
We have to show that $\TW(n)$ is contractible for all $n \geq 0$.
The cases $n=0,1$ are covered by Proposition~\ref{twzeroonepoint}. 
Consider $n \geq 2$.
Use the morphisms $\be_n$, see Definition~\ref{be-n1}, to define a $0$-simplex
$\tau(n) \in  \TW(n)_0=\B_0([0]^n,[0])$  by
\begin{equation}
\label{twpt:tau}
\tau(n) = d_1(\be_n).
\end{equation}
It is depicted by the tree obtained by deleting the strands labelled $1$ from the diagram describing $\be_n$ in Definition~\ref{be-n1}.

For every $k \geq 1$, observe that $(0) \sqcup d^0[k-1]^n=s_1^{k-1}(b_n)$, and define 
\begin{equation}
\label{twpt:Th}
\Theta_{n,k} = s_1^{k-1} (\be_n) \in \B_k([k]^n,(0) \sqcup d^0[k-1]^n).
\end{equation}
We picture $\Theta_{n,k}$ by duplicating $k-1$ times the strand labelled $1$ in 
the figures in Definition~\ref{be-n1} and labelling the resulting strands with $1,\dots,k$.

Augment $\TW$ by $\TW(n)_{-1} = *$. To show that $\TW(n)$ is contractible, it suffices to show that there is a left contraction \cite{dmn:fibrewisecompletion}
\[
s_{-1} \colon \TW(n)_k \to \TW(n)_{k+1}, \qquad k \geq -1
\]
which satisfies 
\begin{equation}\label{twpt:s1}
d_0s_{-1} = \id; \quad d_is_{-1} = s_{-1}d_{i-1}; \quad s_is_{-1} = s_{-1}s_{i-1} \quad \text{ for all } i \geq 0.
\end{equation}
We define $s_{-1} \colon * \to \TW(n)_0$ by $s_{-1}(*) = \tau(n)$.
For any $\vp \in \TW(n)_k $, $k\geq 0$, we let $s_{-1}(\vp)$ be the composition depicted in Figure~\ref{s-1figure}; in formulae:
\[
[k+1]^n \xto{\Theta_{n,k+1}} (0) \sqcup  d^0[k]^n \xto{(0) \sqcup d^0(\vp)} (0) \sqcup d^0[k]= [k+1] \quad \text{(cf. Remark~\ref{cofacial}).}
\]
\begin{figure}[ht]
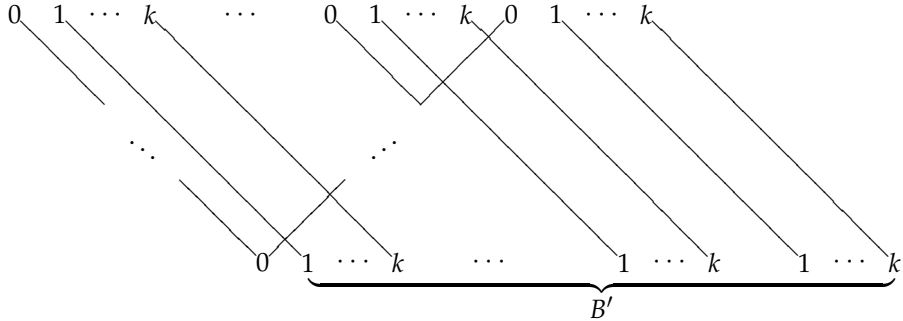

\xy
(0,0)*{0};(12,-12)**\dir{-},
(6,0)*{1};(39,-33)*{1}**\dir{-},
(12,0)*{\cdots},
(18,0)*{k};(51,-33)*{k}**\dir{-},
(30,0)*{\cdots},
(42,0)*{0};(54,-12)**\dir{-},
(48,0)*{1};(81,-33)*{1}**\dir{-},
(54,0)*{\cdots},
(60,0)*{k};(93,-33)*{k}**\dir{-},
(66,0)*{0};(54,-12)**\dir{-},
(72,0)*{1};(105,-33)*{1}**\dir{-},
(78,0)*{\cdots},
(84,0)*{k};(117,-33)*{k}**\dir{-},
(17,-17)*{\ddots},
(22,-22);(33,-33)*{0}**\dir{-};(44,-22)**\dir{-},
(49,-17)*{\iddots},
(45,-33)*{\cdots},
(63,-33)*{\cdots},
(87,-33)*{\cdots},
(111,-33)*{\cdots},
(78,-36)*{\underbrace{\hspace{78mm}}},
(78,-39)*{B'}
\endxy
\caption{The image of $B$ under the simplicial contraction. $B'$ denotes $B$ with all labels increased by $1$.} \label{s-1figure}
\end{figure}
The first identity of \eqref{twpt:s1} holds trivially when $k=-1$. For $k \geq 0$, Lemma \ref{bengames}\eqref{bengames:1}, the simplicial identities, and \eqref{cofacialid} imply
\[
d_0(\Theta_{n,k+1}) = d_0 s_1^k(\be_n) = s_{k+1} \dots s_0^k d_0(\be_n) = \id.
\]
It follows that $d_0s_{-1}=\id$ because for any $\vp \in \TW(n)_k$ where $k \geq 0$ we have
\[
d_0 s_{-1}(\vp) = d_0((0) \sqcup d^0\vp) \circ d_0(\Theta_{n,k+1}) = d_0 d^0(\vp) = \vp.
\]
For the second identity in \eqref{twpt:s1}, note that $k \geq 0$ because $i>0$.
If $k=0$ then $i=1$, and $d_1d^0(\vp)=\id_{()}$ for any $\vp \in \B_0$.
Therefore
\[
d_1s_{-1}(\vp) =   d_1((0) \sqcup d^0(\vp)) \circ d_1(\Theta_{n,1}) =
d_1(\be_n)
= \tau(n) = s_{-1}(*)=s_{-1} d_0(\vp).
\]
If $k\geq 1$, the simplicial identities imply that 
\[
d_i\Theta_{n,k+1}= d_i s_1^k (\be_n) = s_1^{k-1} (\be_n) = \Theta_{n,k}.
\]
Therefore for any $\vp \in \TW(n)_k$, we have
\[
d_i s_{-1}(\vp) = d_i(\Theta_{n,k+1}) \circ d_i((0) \sqcup d^0(\vp)) =
\Theta_{n,k} \circ \{ (0) \sqcup d^0d_{i-1}(\vp)\} = s_{-1}d_{i-1}(\vp).
\]
For the third identity of \eqref{twpt:s1},
we start with the case $i>0$, which forces $k \geq 0$:
\[
s_i (\Theta_{n,k+1}) = s_i s_1^k(\be_n) = s_1^{k+1}(\be_n) = \Theta_{n,k+2}.
\]
It follows that for any $\vp \in \TW(n)_k$,
\[
s_i s_{-1}(\vp) = 
s_i(\Theta_{n,k+1}) \circ s_i((0) \sqcup d^0(\vp)) =
\Theta_{n,k+2} \circ \{ (0) \sqcup d^0s_{i-1}(\vp)\} = s_{-1}s_{i-1}(\vp).
\]
It remains to check that $s_0s_{-1}=s_{-1}s_{-1}$, which is the most subtle identity.
When $k=-1$ we need to prove that $s_{-1}s_{-1}(*) = s_0s_{-1}(*)$.
By the definition of $s_{-1}$ on $\TW(n)_0$ and by Lemma \ref{bengames}\eqref{bengames:4}
\[
s_{-1}s_{-1}(*)=
s_{-1}\tau(n) = 
\{ (0) \sqcup d^0 d_1\be_n\} \circ \be_n = s_0d_1\be_n = s_0\tau(n) = s_0s_{-1}(*).
\]
For $k\geq 0$, consider some $\vp \in \TW(n)_k$. By construction,
\[
s_{-1}(\vp) = \{(0) \sqcup d^0(\vp)\} \circ \Theta_{n,k+1},
\]
and therefore, by \eqref{cofacialid},
\begin{eqnarray*}
s_{0}s_{-1}(\vp) &=& 
    \{(01) \sqcup s_0d^0(\vp)\} \circ s_0(\Theta_{n,k+1}) =
    \{(01) \sqcup d^0d^0(\vp)\} \circ s_0(\Theta_{n,k+1}) 
\\
s_{-1} s_{-1}(\vp) &=&
    \{ (0) \sqcup (1) \sqcup d^0d^0(\vp) \} \circ \{(0) \sqcup d^0(\Theta_{n,k+1}) \} \circ \Theta_{n,k+2}.
\end{eqnarray*}
We are left to show that 
\[
s_0(\Theta_{n,k+1}) = \{(0) \sqcup d^0(\Theta_{n,k+1})\} \circ \Theta_{n,k+2}.
\]
This follows by applying $s_{k+1}\dots s_2$ to Lemma~\ref{bengames}\eqref{bengames:3}
because $\Theta_{n,k+2}=s_1^{k+1}(\be_n)$, and because
\eqref{cofacialid} implies
\begin{eqnarray*}
&& s_0(\Theta_{n,k+1})=s_0 s_1^k (\be_n) = s_2^k s_0(\be_n); \\
&& d^0(\Theta_{n,k+1}) =d^0 s_1^k (\be_n) = s_2^k d^0(\be_n).
\end{eqnarray*}
\end{proof}

%
%
%
%
%
\section{The operad action on completion cosimplicial objects} \label{ringadsection}

Let $(\E,\comp,I)$ be a complete and cocomplete monoidal category. In particular, it is tensored and cotensored over sets by 
\[
E \otimes_\E K = \coprod_K E \qquad \text{ and } \qquad [K,E] = \prod_K E.
\]

\subsection{The simplicial category of cosimplicial objects}

Let $c\E$ denote the category of cosimplicial objects in $\E$. The monoidal structure in $\E$ gives rise to a levelwise monoidal structure $\comp$ in $c\E$.

The category $c\E$ has a simplicial structure \cite[II.2]{goerss-jardine} even if $\E$ does not have one, see \cite[II.1.7ff]{quillen:homalg} or \cite[\S 2.10]{bousfield:cosimplicial-resolutions}.
This structure is called the \emph{external simplicial structure} and is given by
\begin{align*}
[K,X] &= \delta^* [K_\bullet,X^\bullet] & (K \in s\sets,\; X \in c\E)\\
X \otimes_{c\E} K &= (\Lan_\delta X) \otimes_\Delta K & (K \in s\sets,\; X \in c\E)\\
\map_{c\E}(X,Y) &= \Hom_{c\E}(\Lan_\delta X, Y) & (X,\; Y \in c\E),
\end{align*}
where $\de \colon \De \to \De \times \De$ is the diagonal map, $\Lan_\de$ is the left Kan extension along $\de^\op$ and $\otimes_{\De}$ denotes the coend.

If $\De[k]$ denotes the standard $k$-simplex, the usual adjunctions imply
\begin{equation}\label{explicit-k-simplices-mapcE}
\map_{c\E}(X,Y)_k = \Hom_{c\E}(X \otimes_{c\E} \De[k], Y) = \Hom_{c\E}(X , [\De[k],Y]),
\end{equation}
so we think of the $k$-simplices of $\map(X,Y)$ as cosimplicial maps $X \to [\De[k],Y]$.

\begin{prop}
If $(\E,\comp,I)$ is a complete and cocomplete monoidal category then the definitions above make $(c\E,\comp,I)$
a simplicial monoidal category.
\end{prop}
\begin{proof}
Quillen shows in \cite[II.1]{quillen:homalg} that $c\E$ is tensored and cotensored over $s\sets$.

The monoidal structure in $c\E$ also gives rise to
\begin{equation}\label{mon-pairing-ce}
\map_{c\E}(X_1,Y_1) \times \map_{c\E}(X_2,Y_2) \xto{\comp} \map_{c\E}(X_1\comp X_2, Y_1 \comp Y_2)
\end{equation}
as follows.
Given $k$-simplices $\vp \colon X_1 \to [\De[k],Y_1]$ and $\colon X_2 \to [\De[k],Y_2]$ we define 
$\vp \comp \psi \colon X_1\comp X_2 \to [\De[k],Y_1 \comp Y_2]$ in cosimplicial degree $j$ as the composition
$$
X_1^j \comp X_2^j \xto{\prod_{\si} \vp^j(\si) \comp \psi^j(\si)} \prod_{\si \in \De[k]_j} Y_1^j\comp Y_2^j.
$$
One easily checks that $\vp \comp \psi$ is a cosimplicial map and that the assignment $(\vp,\psi) \mapsto \vp \comp \psi$ gives rise to an associative simplicial map \eqref{mon-pairing-ce} which has the constant cosimplicial object $I$ as a unit.

This monoidal product in $c\E$ distributes over the composition. That is, given $k$ simplices $\vp_i \in \map(X_i,Y_i)_k$ and $\psi_i \in \map(Y_i,Z_i)$ ($i=1,\; 2$), we have
\[
(\psi_1 \comp \psi_2) \circ (\vp_1 \comp \psi_2) = (\psi_1 \circ \vp_1) \comp (\psi_2 \circ \vp_2). \qedhere
\]
\end{proof}

Thus, as mentioned in the introduction, every object $X \in c\E$ gives rise to an operad $\opEND(X)$ whose $n$th space is
$\map_{c\E}(X^{\comp n},X)$. Its composition law is obtained from the composition and monoidal structure in $c\E$.

\subsection{Monoids and twist maps}

Before proving Theorem~\ref{twethm}, we show how to derive Corollary~\ref{twcomp} from it.

\begin{proof}[Proof of Corollary~\ref{twcomp}]
We need to see that $t\colon K \comp K \to K \comp K$, given by
\[
t = (K\comp \mu) - \id_{K \comp K} + (\mu \comp K),
\]
is a twist map (Definition~\ref{twistmap}). It is convenient to write $t=d^0s^0 - \id_{K \comp K} +d^1s^0$ where $d^0, d^1$ and $s^0$ are the coface and codegeneracy maps between $K \comp K$ and $K$ in the cosimplicial object $R_K$.

We compute
\begin{eqnarray*}
\mu \circ t = & s^0 \circ (d^0 s^0 - \id +d^1s^0) = s^0 - s^0 - s^0 = s^0  = \mu, \\
t \circ (K\comp \eta) = & (d^0s^0 - \id +d^1s^0) \circ d^0 = d^0 - d^0 + d^1 = d^1 = \eta \comp K\\
(K\comp \mu) \circ (t \comp K) \circ (K \comp t) = &
s^0 \circ  (d^1s^1 - \id_{K \comp K} +d^2s^1) \circ (d^0s^0 - \id_{K \comp K} +d^1s^0) \\
= &d^0s^0s^0-s^1+d^1s^0s^0 = t \circ s^1=t \circ (\mu \comp K).
\end{eqnarray*}
The identities $t \circ (\eta \comp K) = K \comp \eta$ and $t \circ (K \comp \mu) = (\mu \comp K) \circ (K \comp t) \circ (t \comp K)$ follow by symmetry.
\end{proof}

\begin{proof}[Proof of Theorem~\ref{twethm}]
Let $(K,\eta,\mu)$ be a monoid in a complete and cocomplete monoidal category $(\E,\comp,I)$, $R_K \in c\E$ its cobar construction with coface maps $d^i_K$ and codegeneracy maps $s^i_K$. These extend by monoidality to coface and codegeneracy maps on $R_K^{\comp n}$.

It follows from Lemma~\ref{free-monidel-w-relations} and Definition~\ref{def-Bk} of $\B_k$ that there are unique 
monoidal functors $\Phi_k \colon \B_k \to \E$ for all $k \geq 0$ defined by 
\begin{eqnarray*}
&& \Phi_k((a)) = K, \\
&& \Phi_k(u^a) = \eta, \quad \Phi_k(m^a) = \mu, \quad \Phi_k(t^{cd})=t. 
\end{eqnarray*}
The relations (R1) of Definition~\ref{def-Bk} hold due to the fact that $(K,\eta,\mu)$ is a monoid, and the relations 
(R2) hold by Property~\eqref{twmap:2} of the twist map.
On objects, we thus have $\Phi_k([l])=(R_K)^l = K^{\comp (l+1)}$.

In Definition~\ref{def-Bk} it was convenient to define the categories $\B_S$ where $S$ is any finite ordered set.
In this way, given $a \in S$ we defined functors $d_a$ and $s_a$ in Definition~\ref{Bk-simplicial} which render the
categories $\B_k$ a simplicial monoidal category.

We now construct two natural transformations $D^a$ and $S^a$ of functors $\B_S \to \E$. 

\begin{lem}\label{def-Da}
For $a \in S$, there is a natural transformation $D^a \colon \Phi_{d_a(S)} \circ d_a \to \Phi_S$ with
\[
D^a \colon (b) \mapsto 
    \begin{cases}
      I \xto{\eta} K, & \text{if } b=a \\
      K \xto{\id} K, & \text{if } b \neq a,
    \end{cases}
\]
and extended to all objects by monoidality.
Here $d_a(S)=S -\{a\}$ as in Definition~\ref{Bk-simplicial}.
\end{lem}
\begin{proof}
By construction, $\B_k$ defined in \ref{def-Bk} is monoidally generated by the objects $(b)$ and the morphisms $u^b, m^b$ and $t^{cd}$.
To prove that $D^a$ is a natural transformation, it suffices by Lemma~\ref{free-monidel-w-relations} to show that
for all $a,b,c,d \in S$,
\begin{eqnarray*}
&& \Phi_S(u^b)\circ D^a(()) = D^a((b)) \circ \Phi_{d_a(S)}( d_a(u^b)), \\
&& \Phi_S(m^b) \circ D^a((bb)) = D^a((b)) \circ \Phi_{d_a(S)}(d_a(m^b)), \\ 
&& \Phi_{S}(t^{cd}) \circ D^a((cd)) = D^a((dc)) \circ \Phi_{d_a(S)} (d_a(t^{cd})).
\end{eqnarray*}
The first equality hold if $a=b$ because $d_a(u^a)=\id$.
If $a \neq b$ then $D^a((b))=\id=D^a(())$ and the equality holds again since $d^a(u^b)=u^b$.

The second equality holds if $a\neq b$ because $D^a((bb))=\id=D^a((b))$ and $d^a(m^b)=m^b$.
If $a =b$ then $d_a(m^a)=\id_{()}$ so by the relations between $\eta$ and $\mu$,
$$
\Phi_S(m^a) \circ D^a((aa)) = \mu \circ (\eta\comp \eta) = \eta = D^a((a)).
$$

The third equality holds when $a \neq c,\; d$ because $D^a((cd)) = \id_{K\comp K} = D^a((dc))$ and $d_a(t^{cd})=t^{cd}$.
If $a=c$ then $d_a(t^{cd})=\id$ (cf. Definition~\ref{def-Bk}) so 
\[
\Phi_S(t^{cd}) \circ D^a((cd)) =
t \circ (\mu \comp K) = K \comp \eta = D^a((dc)) \circ \Phi_{d_aS} (d_a(t^{cd})).
\]
A similar calculation applies when $a =d$.
\end{proof}

\begin{lem}\label{def-Sa}
For $a \in S$, there is a natural transformation $S^a \colon \Phi_{s_a(S)} \circ s_a \to \Phi_S$ with
\[
S^a \colon (b) \mapsto \begin{cases}
      K \comp K \xto{\mu} K, & \text{if } b=a \\
      K \xto{\id} K, & \text{if } b \neq a,
    \end{cases}
\]
and extended to all objects by monoidality.
Here $s_a(S)=S -\{a\} \sqcup \{a',a''\}$ as in Definition~\ref{Bk-simplicial}.
\end{lem}
\begin{proof}
As in the proof of Lemma ~\ref{def-Da}, in order to prove the naturality of $S^a$ we have to show that
\begin{eqnarray*}
\Phi_S (u^b) \circ S^a() &=& S^a(b) \circ \Phi_{s_aS} (s_a(u^b)), \\
\Phi_S (m^b) \circ S^a(bb) &=& S^a(b) \circ \Phi_{s_aS} (s_a(m^b)), \\
\Phi_S (t^{cd}) \circ S^a(cd) &=& S^a(dc) \circ \Phi_{s_aS} (s_a(t^{cd})).
\end{eqnarray*}
The first and second equalities are easy to prove when $a \neq b$ because $S^a(bb)=\id=S^a(b)$ and $s^a(m^b)=m^b$.
When $a=b$ we see from Definition~\ref{def-Bk} that
\[
\Phi_S(u^b) \circ S^a() = u = \mu \circ (\eta \comp \eta)  = S^a(a) \circ \Phi_{s_aS} (s_a(u^a)).
\]
For the second equality ($a=b$), we use the associativity of $\mu$ and Definitions~\ref{twistmap}\eqref{twmap:comm} and \ref{Bk-simplicial} to deduce that 
\begin{align*}
& S^a(a) \circ \Phi_{s_aS}(s^a(m^a)) = \mu \circ \Phi_{s_aS}(m^{a'} \sqcup m^{a''}) \circ \Phi_{s_aS}((a')\sqcup t^{a''a'} \sqcup (a''))\\
= & \mu \circ (\mu \comp \mu) \circ (K \comp t \comp K) = \mu \circ (\mu \comp \mu) = \Phi_S(m^a) \circ S^a(aa).
\end{align*}
The third equality is again straightforward when $a\neq c,\; d$.
If $a=c$ then by Definitions~\ref{def-Bk}, \ref{Bk-simplicial} and \ref{twistmap}\eqref{twmap:3},
\begin{align*}
& S^a(da) \circ \Phi_{s_aS}( s_a(t^{ad})) = 
(K \comp \mu) \circ \Phi_{s_aS}(t^{a'd}\sqcup (a'')) \circ \Phi_{s_aS}((a') \sqcup t^{a''d})\\
= & (K \comp \mu) \circ (t \comp K) \circ (K \comp t) = t \circ (\mu \comp K) =  
\Phi_S(t^{ad}) \circ S^a(ad).
\end{align*}
A similar calculation applies when $a=d$.
\end{proof}

When $S=[k]$, we observe that $\Phi_S([k]) = K^{\comp (k+1)} = R_K^k$ and that for any $i \in [k]$ we have $D^i([k]) = d^i_{R_K}$, $S^i([k]) = s^i_{R_K}$, and consequently $D^i([k]^n) = d^i_{R_K^{\comp n}}$ and $S^i([k]^n) = s^i_{R_K^{\comp n}}$.

\begin{defn}
For integers $n,m \geq 0$ let $\TW(n,m)$ denote the simplicial set whose set of $k$-simplices is $\B_k([k]^n,[k]^m)$. In particular, $\TW(n,1) = \TW(n)$ (cf. Definition~\ref{def-tw}).
\end{defn}

Here we use the fact that $d_i \colon \B_{k} \to \B_{k-1}$ and $s_i \colon  \B_{k} \to \B_{k+1}$ are monoidal functors which carry $[k]$ to $[k-1]$ and $[k+1]$ respectively, see Proposition \ref{bk-simplcl}. 
Composition of morphisms and the monoidal operation $\sqcup$ in $\B_k$ give rise to ``composition'' and ``monoidal'' 
products of simplicial sets
\begin{eqnarray}\label{compmonpair}
&& \TW(n,m) \times \TW(m,l) \xto{\circ} \TW(n,l), \\
\nonumber
&&\TW(n,m) \times \TW(n',m') \xto{\sqcup} \TW(n+n',m+m').
\end{eqnarray}

Write $X(n) = R_K^{\comp n}$ and $E(m,n) = \map_{c\E}(X(m),X(n))$. Thus $\Phi_k([k]^n)=X(n)^k$. By  Proposition~\ref{bk-simplcl}, $\B_\bullet$ is a simplicial monoidal category and thus a morphism $\al \in \De([i],[j])$ gives rise to a functor $\al^* \colon \B_j \to \B_i$.

For $f \in \TW(n,m)_k$ and $j\geq 0$, we define a morphism $X(n)^j \to [\De[k], X(m)]^j$ by
$$
X(n)^j \xto{ \quad \prod_{\al} \Phi_j(\al^*(f)) \quad } \prod_{\al \in \De[k]_j} X(m)^j.
$$
These maps assemble to a cosimplicial map 
\[
\rho(n,m)_k(f) \colon X(n) \to [\De[k],X(m)],
\] 
because by Lemma~\ref{def-Da}, for every injective $\epsilon^i \colon [j] \to [j+1]$ in $\De$ and every 
$\al \in \De[k]_{j+1}$,
\begin{multline*}
d^i_{X(m)} \circ \Phi_j((\al \circ \epsilon^i)^*(f)) = D^i([k]^m) \circ \Phi_{d_i[j+1]}(d_i(\al^*(f))) = \\
\Phi_{j+1}(\al^*(f)) \circ D^i([k]^n) = \Phi_{j+1}(\al^*(f)) \circ d^i_{X(n)}.
\end{multline*}
Similarly, by Lemma~\ref{def-Sa}, for any  surjective $\eta^i \colon [j] \to [j-1]$ and every $\al \in \De[k]_{j-1}$,
\begin{multline*}
s^i_{X(m)} \circ \Phi_j((\al \circ \eta^i)^*(f)) = S^i([k]^m) \circ \Phi_{s_i([j-1])}(s_i(\al^*(f)) = \\
 \Phi_{j-1}(\al^*(f)) \circ S^i([k]^n) = \Phi_{j-1}(\al^*(f)) \circ s^i_{X(n)}.
\end{multline*}
Thus, $\rho(n,m)_k(f) \in E(n,m)_k$, see \eqref{explicit-k-simplices-mapcE}. They assemble to form a simplicial map
\[
\rho(n,m) \colon \TW(n,m) \to E(n,m).
\]
because for any $\ga \in \De([\ell],[k])$, $\al \in \De[\ell]_j$, and $f \in \TW(n,m)_k$, we have
$$
\Phi_j((\ga \circ \al)^*(f)) = \Phi_j(\al^*(\ga^*(f)),
$$
so $\rho(n,m)_\ell(\ga^*(f))=\ga^*\rho(n,m)_k(f)$.

The following diagrams of spaces commute:
$$
\xymatrix{
\TW(n,m) \times \TW(m,l) \ar[d]_{\rho(n,m) \times \rho(m,l)} \ar[r]^-{\circ} &
\TW(n,l) \ar[d]^{\rho(n,l)}
\\
E(n,m) \times E(m,l) \ar[r]_-{\circ} &
E(n,l)
}
$$
$$
\xymatrix{
\TW(n,m) \times \TW(n',m') \ar[d]_{\rho(n,m) \times \rho(n',m')} \ar[r]^{\sqcup} &
\TW(n+n',m+m') \ar[d]^{\rho(n+n',m+m')} 
\\
E(n,m) \times E(n',m') \ar[r]_{\comp}^{\eqref{mon-pairing-ce}} &
E(n+n',m+m').
}
$$
This follows by verification in every simplicial degree $k$ using the fact that the functors $\Phi_k$
are monoidal.
For example, given $f \in \TW(n,m)_k$ and $g \in \TW(m,l)_l$, the cosimplicial map $\rho(m,l)_k(g) \circ \rho(n,m)_k(f)$ has the form
$$
X(n)^j \xto{\prod_{\al} \Phi_j(\al^*(g)) \circ \Phi(\al^*(f))} \prod_{\al \in \De[k]_j} X(l)^j.
$$
Now $\Phi_j(\al^*(g)) \circ \Phi_j(\al^*(f)) = \Phi_j(\al^*(g \circ f))$; which assemble to $\rho(n,l)_k(g \circ f)$.
A similar argument works for the second diagram.

Since the composition law in $\TW$ (Definition~\ref{def-tw}) is derived from \eqref{compmonpair}, the commutativity of the diagrams above implies that the maps
\[
\rho(n,1) \colon \TW(n) \to \map_{c\E}(X(n),X(1))
\]
form a morphism of operads $\TW \to \opEND(R_K)$. 
This completes the proof of Theorem~\ref{twethm}.
\end{proof}

\appendix

\section{The free monoidal construction}

In this appendix we will prove Lemma~\ref{free-monoidal}. The result is an exercise in free constructions, and probably known to anybody who has worked with small monoidal categories, but we could not find it in the literature.

\begin{lem}
The functor $U_{vm}\colon \{\text{monoidal graphs}\} \to \{\text{vertex-monoidal graphs}\}$ has a left adjoint.
\end{lem}
\begin{proof}
The adjoint is given by $F_{vm}(\A \rightrightarrows \O) = (\A_\sqcup \rightrightarrows \O)$, where $\A_\sqcup$ is the free monoid on $\A$ and the structure maps are extended uniquely to monoidal maps. The adjointness is obvious.
\end{proof}

\begin{lem}
The functor $U_m\colon \{\text{unital monoidal graphs}\} \to \{\text{monoidal graphs}\}$ has a left adjoint.
\end{lem}
\begin{proof}
Let $\A \rightrightarrows \O$ be a monoidal graph. For every $\sigma \in \O$ we introduce a symbol $\iota(\sigma) \notin \A$ and let $\tilde \Delta$ denote $\A \sqcup \{ \iota(\sigma) \}_{\sigma \in \O}$.
Extend $\partial_0,\partial_1\colon \A \to \O$ to $\tilde \Delta$ by setting $\partial_0(\iota(\sigma))=\partial_1(\iota(\sigma))=\sigma$.

This is not quite a unital monoidal graph because $\iota\colon \O \to \tilde \Delta$ is not monoidal. We enforce this by setting $\Delta = \tilde \Delta/\sim$, where $\sim$ is the equivalence relation generated by
\[
\iota(1_\O) \sim 1_\A \quad \text{and} \quad \iota(\sigma \sqcup \tau) \sim \iota(\sigma) \sqcup \iota(\tau).
\]
Since $\partial_i(a) = \partial_i(a')$ if $a \sim a' \in \tilde\Delta$, they are well-defined on$\Delta$, and we obtain a unital monoidal graph $(\Delta,\O)$. Again, adjointness is easily checked.
\end{proof}

\begin{lem}
The functor $\{\text{small monoidal categories}\} \xto{\, U_{um} \,} \{\text{unital monoidal graphs}\}$ has a left adjoint.
\end{lem}
\begin{proof}
Let $\G = (\partial_0,\;\partial_1\colon \A \overset{\curvearrowleft}{\rightrightarrows} \O\!\iota)$ be a unital monoidal graph. Let $\tilde{\M}$ denote $C(\G)$, the free category generated by the graph $\G$. However, $\tilde{\M}$ is not monoidal yet.

Let $\sim$ be the smallest equivalence relation on the morphism set of $\tilde{\M}$ with
\begin{enumerate}\renewcommand{\theenumi}{\roman{enumi}}
\item \label{compositionmonoidalrel} $(f \sqcup g) \circ (f' \sqcup g') \sim (f \circ f') \sqcup (g \circ g')$ for all composable arrows $f,f'$ and $g,g'$ in $\tilde{\M}$,
\item \label{closuremonoidalrel} $\sim$ is closed under $\sqcup$, i.e. if $f \sim f'$ and $g \sim g'$ then $(f \sqcup g) \sim (f' \sqcup g')$.
\end{enumerate}
Again, the equivalence relation $\sim$ has the property that if $f \sim f'$ then their domains and codomains are equal, and thus we obtain a quotient category $\M = \tilde{M}/\sim$ \cite[\S II.8, Proposition 1]{maclane:cwm}.
The monoidal structure $\sqcup$ descends from $\tilde{M}$ thanks to \eqref{closuremonoidalrel}.
It renders $\M$ a monoidal category by virtue of relation \eqref{compositionmonoidalrel}.
\end{proof}


\bibliographystyle{alpha}
\bibliography{bibliography}

\end{document}